\documentclass[12pt]{amsart}
\usepackage{hyperref}
\date{23rd September 2014}
\DeclareMathOperator\End{End}
\newcommand\ii{\mathbf i}
\newcommand\jj{\mathbf j}
\newcommand\kk{\mathbf k}
\newcommand\n{\mathbf n}

\newtheorem{theorem}{Theorem}

\theoremstyle{remark}
\newtheorem{remark}{Remark}
\newtheorem{example}{Example}

\newcommand{\mat}[4]{\begin{pmatrix}#1&#2\\#3&#4\end{pmatrix}}
\newcommand{\smat}[4]{\left(\begin{smallmatrix}#1&#2\\#3&#4\end{smallmatrix}\right)}

\title{Centre of the Schur Algebra}
\author{T. Geetha}\address{Chennai Mathematical Institute, Siruseri}
\author{Amritanshu Prasad}\address{The Institute of Mathematical Sciences, Chennai}
\begin{document}
\maketitle
\begin{abstract}
  We describe a basis of the centre of the Schur algebra which comes from conjugacy classes in the symmetric group via Schur-Weyl duality.
  We give a combinatorial description of expansions of these basis elements in terms of the basis originally used by Schur.
  The primitive central idempotents of the Schur algebra can be written down using this basis and the character table of the symmetric group.
  Along the way we prove a result on the non-singularity of the submatrix of the character table matrix of a symmetric group obtained by taking rows and columns indexed by partitions with at most $n$ parts for any $n$.
\end{abstract}
\section{The Schur Algebra}
\label{sec:schur-alg}

Let $n$ and $d$ be two positive integers.
Let $\n$ denote the set $\{1,\dotsc,n\}$ of the first $n$ positive integers.
Let $I(n,d) = \n^d$, the $d$-fold Cartesian power of $\n$.
The symmetric group $S_d$ on $d$ symbols acts on an element $\ii = (i_1,\dotsc,i_d)$ of $I(n,d)$ by permuting the $d$-components of an element of $I(n,d)$:
\begin{equation*}
  w\cdot (i_1,\dotsc,i_d) = (i_{w(1)}, \dotsc,i_{w(d)}) \text{ for each }w\in S_d. 
\end{equation*}
Let $S_d$ acts on $I(n,d)^2$ by the diagonal action:
\begin{equation*}
  w\cdot (\ii,\jj) = (w\cdot \ii, w\cdot \jj) \text{ for all } w\in S_d, \ii,\jj \in I(n,d).
\end{equation*}
Let $M(n,d)$ denote the set of $n\times n$ matrices whose entries are non-negative integers which add up to $d$.
Define a function $D:I(n,d)^2\to M(n,d)$ as follows:
\begin{equation*}
  D(\ii,\jj)_{i,j} = \#\{1\leq k\leq d\mid (i_k,j_k) = (i,j)\} \text{ for all } i,j\in \n.
\end{equation*}
The following is straightforward:
\begin{theorem}
  \label{thorem:bijection}
  The function $D:I(n,d)^2\to M(n,d)$ defined above gives rise to a bijection from the set of $S_d$-orbits in $I(n,d)^2$ onto $M(n,d)$.
\end{theorem}
According to Green \cite[Section~2.3]{Green}, the Schur algebra $S_K(n,d)$ has a basis indexed by $S_d$-orbits in $I(n,d)^2$.
Writing $\xi_D$ for the basis element corresponding to the orbit consisting of pairs $(\ii, \jj)$ such that $D(\ii,\jj) = D$ for each $D\in M(n,d)$, Green\footnote{Instead of $\xi_D$, Green uses the notation $\xi_{\ii,\jj}$ where $(\ii,\jj)\in I(n,d)^2$ is any pair such that $D(\ii,\jj)=D$.} describes the structure constants of the Schur algebra as follows:
let $D$, $D'$ and $D''$ be three matrices in $M(n,d)$.
Choose a pair $(\ii,\jj)\in I(n,d)^2$ such that $D(\ii,\jj) = D$.
Then,
\begin{equation*}
  \xi_{D'} \xi_{D''} = \sum_\kk c^D_{D',D''} \xi_D,
\end{equation*}
where
\begin{equation*}
  c^D_{D',D''} = |\{\kk\in I(n,d)\mid D(\ii,\kk) = D' \text{ and } D(\kk,\jj) = D''\}|.
\end{equation*}
For more on these structure constants, see \cite{structure}.

\section{Permutation Representations}
\label{sec:perm-repr}

Let $G$ be a finite group and $X$ be a finite $G$-set.
Let $K[X]$ denote the set of all $K$-valued functions on $X$.
We may think of $K[X]$ as a representation of $G$ using the action:
\begin{equation*}
  \rho_X(g)f(x) = f(g^{-1}\cdot x).
\end{equation*}

Given a function $k:X\times X\to K$, define $T_k\in \End_K K[X]$ by
\begin{equation}
  \label{eq:intop}
  T_k f(x) = \sum_{y\in X} k(x,y)f(y).
\end{equation}
The operator $T_k$ is called the integral operator associated to the kernel $k$.

Now any linear endomorphism $T:K[X]\to K[X]$ can be expressed as an \emph{integral operator} with respect to a unique integral kernel $k:X\times X\to K$; in other words, there exists a unique function $k:X\times X\to K$ such that $T=T_k$.
To find $k$, let $\delta_x\in K[X]$ denote the function that is $1$ at $x$ and $0$ everywhere else.
Then if 
\begin{equation*}
  T\delta_y = \sum_{x\in X} T_{xy}\delta_x,
\end{equation*}
then setting $k(x,y) = T_{xy}$ for all $x,y\in X$, one easily verifies that $T = T_k$.

Let
\begin{equation*}
  K[X\times X]^G = \{k:X\times X\to K\mid  k(g\cdot x, g\cdot y) = k(x,y) \text{ for }g\in G, x,y\in X\}.
\end{equation*}
Among the linear operators on $K[X]$, the $G$-endomorphisms are characterized by
\begin{equation*}
  \End_G K[X] = \{T_k\mid k\in K[X\times X]^G\}. 
\end{equation*}

\section{Tensor Space as a Permutation Representation}
\label{sec:tensor-space-as}

Since $I(n,d)$ comes with an $S_d$-action (see Section~\ref{sec:schur-alg}),
we may view $K[I(n,d)]$ as a permutation representation of $S_d$.

Take $V=K^n$ with coordinate vectors $e_1,\dotsc,e_n$.
The symmetric group $S_d$ acts on $V^{\otimes d}$ by:
\begin{equation*}
  w\cdot(v_1\otimes \dotsb\otimes v_d) = v_{w(1)}\otimes \dotsb \otimes v_{w(d)}.
\end{equation*}
For each $\ii=(i_1,\dotsc,i_d)\in I(n,d)$, consider the vector in $V^{\otimes d}$
\begin{equation*}
  e_\ii = e_{i_1}\otimes \dotsb \otimes e_{i_d}.
\end{equation*}
Then
\begin{equation*}
  \{e_\ii\mid \ii \in I(n,d)\}
\end{equation*}
is a basis of $V^{\otimes d}$.
The following theorem is straightforward:
\begin{theorem}
  \label{theorem:tensor-as-perm}
  The map $\delta_\ii\mapsto e_\ii$ defines an isomorphism of $S_d$-representations $K[I(n,d)]\to V^{\otimes d}$.
\end{theorem}
In particular, for each $w\in S_d$,
\begin{equation*}
  w\cdot e_\jj = e_{w\cdot \jj}.
\end{equation*}
When one views $V^{\otimes d}$ as $K[I(n,d)]$ under the isomorphism of Theorem~\ref{theorem:tensor-as-perm}, we may write:
\begin{equation*}
  w\cdot \delta_\jj = \delta_{w\cdot \jj}.
\end{equation*}
Therefore, by the discussion in Section~\ref{sec:perm-repr}, the operator by which $w$ acts on $K[I(n,d)]$ is $T_{k_w}$, where $k_w:I(n,d)\times I(n,d)\to K$ is the function
\begin{equation}
  \label{eq:kw}
  k_w(\ii,\jj) =
  \begin{cases}
    1 & \text{if } \ii = w\cdot \jj,\\
    0 & \text{otherwise.}
  \end{cases}
\end{equation}
On the other hand, if $D$ is an $n\times n$ integer matrix with sum $d$, then the element $\xi_D\in S_K(n,d)$ acts on $V^{\otimes d}$ by
\begin{equation*}
  \xi_D e_\jj = \sum_{\{\ii \mid D(\ii,\jj)=D\}} e_\ii.
\end{equation*}
When $V^{\otimes d}$ is viewed as $K[I(n,d)]$, then we may write
\begin{equation*}
  \xi_D \delta_\jj = \sum_{\{\ii \mid D(\ii,\jj)=D\}} \delta_\ii.
\end{equation*}
Thus, $\xi_D = T_{k_D}$, where $k_D:I(n,d)\times I(n,d)\to K$ is given by
\begin{equation}
  \label{eq:kD}
  k_D(\ii,\jj) =
  \begin{cases}
    1 & \text{if } D(\ii,\jj) = D,\\
    0 & \text{otherwise.}
  \end{cases}
\end{equation}
For each partition $\lambda$ of $d$, let
\begin{equation}
  \label{eq:cla-defn}
  c_\lambda = \sum_{w\text{ has cycle type }\lambda} 1_w,
\end{equation}
where $1_w$ is the basis element of $K[S_d]$ corresponding to $w\in S_d$.
Since cycle types determine the conjugacy classes of $S_d$, the set
\begin{equation*}
  \{c_\lambda \mid \lambda \text{ is a partition of }d\}
\end{equation*}
is a basis of the centre of $K[S_d]$.
Note that the centre of $K[S_d]$ maps onto the centre of $S_K(n,d)$.
For an element $z$ in the centre of $K[S_d]$, let $\bar z$ denote its image in the centre of $S_K(n,d)$.
Consider the coefficients $c_{\lambda,D} \in K$ of the expansion of $\bar c_\lambda$ in terms of Schur's basis for $S_K(n,d)$:
\begin{equation}
  \label{eq:cla-expansion}
  \bar c_\lambda = \sum_{D\in M(n,d)} c_{\lambda,D}\xi_D.
\end{equation}

Now the action of $c_\lambda$ on $K[I(n,d)]$ is given by the kernel
\begin{equation}
  \label{eq:kla-kw}
  k_\lambda(\ii,\jj) = \sum_{w \text{ has cycle type }\lambda} k_w(\ii,\jj).
\end{equation}
On the other hand, by equating the kernels of the operators on both  sides of (\ref{eq:cla-expansion}), we get:
\begin{equation}
  \label{eq:kla-kD}
  k_\lambda(\ii,\jj) = \sum_{D\in M(n,d)} c_{\lambda,D}k_D(\ii,\jj).
\end{equation}
Let $(\ii,\jj)\in I(n,d)^2$ be any pair with $D(\ii,\jj) = D$.
Then evaluating the right hand sides of the equations (\ref{eq:kla-kw}) and (\ref{eq:kla-kD}) at $(\ii,\jj)$ and equating them gives 
\begin{equation*}
  c_{\lambda,D} = \#\{w\in S_d\mid w\text{ has cycle type }\lambda \text{ and } w\cdot \jj =\ii\}.
\end{equation*}
Thus we have proved the following theorem:
\begin{theorem}
  \label{theorem:claD}
  For each partition $\lambda$ of $d$, let $c_\lambda$ be the element of the centre of $K[S_d]$ corresponding to permutations of cycle type $\lambda$, as in Eq.~(\ref{eq:cla-defn}).
  Then the image $\bar c_\lambda$  of $c_\lambda$ in $\End(V^{\otimes d})$ of $c_\lambda$ has expansion
  \begin{equation*}
    \bar c_\lambda = \sum_{D\in M(n,d)} c_{\lambda,D}\xi_D,
  \end{equation*}
  where, for any pair $(\ii,\jj)\in I(n,d)^2$ with $D(\ii,\jj) = D$,
  \begin{equation*}
    c_{\lambda,D} = \#\{w\in S_d\mid w\text{ has cycle type } \lambda \text{ and } w\cdot \jj = \ii\}.
  \end{equation*}
\end{theorem}
\begin{remark}
  In order to use Theorem~\ref{theorem:claD}, we need an algorithm to construct a pair $(\ii,\jj)$ with $D(\ii,\jj) = D$ for a given $D\in M(n,d)$.
  For this, let
  \begin{equation*}
    \mu_i = \sum_{j=1}^n D_{ij} \text{ and } \nu_j = \sum_{i=1}^n D_{ij}
  \end{equation*}
  be the row and column sums of $D$, respectively.
  Let $\mu(D) = (\mu_1,\dotsc,\mu_n)$ and $\nu(D) = (\nu_1,\dotsc,\nu_n)$.
  These are weak compositions of $d$ with $n$ parts\footnote{By definition, a weak composition of $d$ with $n$ parts is a tuple $(\lambda_1,\dotsc,\lambda_d)$ of non-negative integers such that $\lambda_1+\dotsb+\lambda_d = d$.}.
  Take $\jj$ to be the vector:
  \begin{equation*}
    \jj = (\underbrace{1,\dotsc,1}_{\nu_1 \text{ times}},\underbrace{2,\dotsc,2}_{\nu_2 \text{ times}},\dotsc,\underbrace{n,\dotsc,n}_{\nu_n \text{ times}}).
  \end{equation*}
  Take $\ii$ to be the vector whose coordinates corresponding to the $\nu_j$ instances of $j$ in $\jj$ are given by
  \begin{equation*}
    \dotsc,\underbrace{1,\dotsc,1}_{D_{1j} \text{ times}},\underbrace{2,\dotsc,2}_{D_{2j} \text{ times}},\dotsc,\underbrace{n,\dotsc,n}_{D_{nj} \text{ times}},\dotsc
  \end{equation*}
  Then $D(\ii,\jj) = D$.

  For example, if
  \begin{equation*}
    D =
    \begin{pmatrix}
      1 & 1 & 0 & 0\\
      0 & 0 & 0 & 0\\
      0 & 1 & 0 & 2\\
      0 & 0 & 0 & 0
    \end{pmatrix}
    ,
  \end{equation*}
  then $\mu(D) = (2, 0, 3, 0)$ and $\nu(D) = (1, 2, 0, 2)$.
  We have
  \begin{equation*}
    \jj = (1, 2, 2, 4, 4)
  \end{equation*}
  and
  \begin{equation*}
    \ii = (1, 1, 3, 3, 3).
  \end{equation*}
\end{remark}
In the example above, there does not exist and $w\in S_d$ such that $\ii = w\cdot \jj$.
In fact, for general $\ii$ and $\jj$ such permutations do not exist unless each of the integers $1,2,\dotsc,n$ occurs in $\ii$ the same number of times that it occurs in $\jj$.
In other words, the weak compositions $\mu(D)$ and $\nu(D)$ associated to the matrix $D$ are the same.
\begin{theorem}
  Each central element $\xi\in S_K(n,d)$ lies in the span of
  \begin{equation*}
    \{\xi_D\in M(n,d)\mid \mu(D) = \nu(D)\}.
  \end{equation*}
\end{theorem}

\section{Central idempotents and the action of $K[S_n]$}
\label{sec:central-idempotents}
For any finite dimensional split semisimple unital $K$-algebra $A$, let
\begin{equation*}
  1 = \epsilon_1 + \epsilon_2 + \dotsb + \epsilon_r
\end{equation*}
be the decomposition of the multiplicative unit into a sum of primitive central idempotents.
The isomorphism classes of simple $A$-modules can be enumerated as $V_1,\dotsc, V_r$ where $\epsilon_i$ acts on $V_i$ by the identity map, and on $V_j$ by $0$ for all $j\neq i$
.
Now if $W$ is any finite  dimensional $A$-module whatsoever, whose decomposition into simple modules with multiplicities is given by
\begin{equation}
  \label{eq:decomp-simple-with-mults}
  W = V_1^{\oplus m_1}\oplus\dotsb \oplus V_r^{\oplus m_r},
\end{equation}
then, since $\epsilon_i$ acts on $V_j$ by $0$ unless $i=j$, those $\epsilon_i$'s for which $m_i=0$ lie in the kernel of the map $A\to \End_KW$.
On the other hand, the $\epsilon_i$'s for which $m_i>0$ have non-zero image in $\End_KW$.
Also, since the image of any element in the centre of $A$ lies in $\End_AW$, Schur's lemma implies that such an image lies in the span of the images of the $\epsilon_i$'s for which $m_i>0$.
Thus we have
\begin{theorem}
  \label{theorem:image-of-alg}
  Let $A$ be a split semisimple algebra over a field $K$.
  Let $W$ be an $A$-module whose decomposition into a sum of simple modules with multiplicities is given by~\eqref{eq:decomp-simple-with-mults}.
  Then the image of the centre of $A$ in $\End_KW$ is spanned by the images of
  \begin{equation*}
    \{\epsilon_i\mid 1\leq i\leq r,\: m_i>0\}.
  \end{equation*}
  Also, for each $i$ such that $m_i=0$, $\epsilon_i$ lies in the kernel of the map $A\to \End_K W$.
\end{theorem}

Now consider the case where $A = K[S_d]$ and $W = K[I(n,d)]$.
The simple modules of $K[S_d]$ are indexed by partitions of $d$.
For each partition $\lambda$ of $n$ write $\epsilon_\lambda$ for the corresponding primitive central idempotent in $K[S_d]$.

For each weak composition $\lambda=(\lambda_1,\dotsc,\lambda_n)$ of $d$ with $n$ parts, set
\begin{equation*}
  X_\lambda = \{\ii\in I(n,d)\mid i_k = i \text{ for exactly $\lambda_i$ values of $k$}\}.
\end{equation*}
Then
\begin{equation*}
  I(n,d) = \coprod_\lambda X_\lambda
\end{equation*}
is the decomposition of $I(n,d)$ into orbits for the action of $S_d$.
Thus, as a representation of $S_d$,
\begin{equation*}
  K[I(n,d)] = \bigoplus_\lambda K[X_\lambda],
\end{equation*}
the sum being over the weak compositions of $d$ with $n$ parts.

\begin{example}
  Take $n=2$ and $d=3$. There are four weak compositions of $3$ with two parts, namely,
  \begin{equation*}
    (3, 0), (0,3), (2,1) \text{ and } (1,2).
  \end{equation*}
  Thus $I(2,3)$ has four $S_3$-orbits.
  The orbit corresponding to each weak composition is given in the table below:
  \begin{equation*}
    \begin{array}{cc}
      \hline
      \lambda & X_\lambda\\
      \hline
      (3,0) & \{(1,1,1)\}\\
      (0,3) & \{(2,2,2)\}\\
      (2,1) & \{(1,1,2), (1,2,1), (2,1,1)\}\\
      (1,2) & \{(2, 2, 1), (2, 1, 2), (1, 2, 2)\}\\
      \hline
    \end{array}
  \end{equation*}
\end{example}

If $\tilde\lambda$ is obtained from $\lambda$ by permuting its coordinates, then $X_\lambda$ and $X_{\tilde\lambda}$ are isomorphic as $S_d$-sets (there is an obvious bijection from $X_\lambda$ onto $X_{\tilde\lambda}$ which respects the $S_d$ actions on these sets).
It follows that the permutation representations $K[X_\lambda]$ and $K[X_{\tilde\lambda}]$ are isomorphic.

Thus, a representation of $S_d$ occurs in $K[I(n,d)]$ if and only if it occurs in $K[X_\lambda]$ for some partition $\lambda$ of $d$ with at most $n$ parts (a partition with fewer than $n$ parts can be thought of as a weak composition with $n$ parts by appending $0$'s at the end).
The representation $K[X_\lambda]$ is nothing but the representation of $S_n$ induced from the trivial representation of the Young subgroup $S_{\lambda_1}\times \dotsb \times S_{\lambda_n}$, which by Young's rule \cite[Theorem~3.3.1]{rtcv} has decomposition
\begin{equation*}
  K[X_\lambda] = \bigoplus_{\mu\leq \lambda} V_\mu^{\oplus K_{\mu\lambda}}.
\end{equation*}
Here $\mu\leq \lambda$ signifies that $\mu$ is less than or equal to $\lambda$ in the reverse dominance order, which is the same as saying that
\begin{equation*}
  \mu_1+\dotsb + \mu_k \geq \lambda_1 + \dotsb + \lambda_k
\end{equation*}
for each $k$.
Also, it is known that $K_{\mu\lambda}>0$ whenever $\mu \leq \lambda$ \cite[Theorem~3.1.12]{rtcv}.

We may conclude that $V_\mu$ occurs in $K[I(n,d)]$ if and only if $\mu\leq \lambda$ for some partition $\lambda$ with at most $m$ parts.
It is also well-known that if $\mu\geq \lambda$, then the number of non-zero parts in $\mu$ is at most the number of non-zero parts in $\lambda$ \cite[Exercise~3.1.11]{rtcv}.
Thus $V_\mu$ occurs in $K[I(n,d)]$ if and only if $\mu$ has at most $n$ parts.

Combining this with Theorem~\ref{theorem:image-of-alg} (and remembering that $K[I(n,d)]$ is nothing but the tensor space $V^{\otimes d}$, we obtain the following result:
\begin{theorem}
  \label{theorem:basis-of-image}
  For positive integers $n$ and $d$, the image of the centre of $K[S_d]$ in $\End_KV^{\otimes d}$ is spanned by
  \begin{equation*}
    \{\bar \epsilon_\lambda\mid \lambda\in P(n,d)\},
  \end{equation*}
  where $P(n,d)$ denotes the set of partitions of $d$ with at most $n$ parts.
  Moreover, if $\lambda$ is a partition of $d$ with more than $n$ parts, then $\bar \epsilon_\lambda =0$.
\end{theorem}

\section{A basis for the centre of the Schur algebra}
\label{sec:basis-centre-schur}
We know that the functions $c_\lambda$ defined in (\ref{eq:cla-defn}) form a basis for the centre of $K[S_d]$.
For its image in $\End_K V^{\otimes d}$, we have the following theorem:
\begin{theorem}
  \label{theorem:class-basis-of-center}
  The set
  \begin{equation}
    \label{eq:basis-of-class-fns}
    \{\bar c_\lambda\mid \lambda\in P(n,d)\}
  \end{equation}
  is a basis for the centre of the image of $K[S_d]$ in $\End_K V^{\otimes d}$.
\end{theorem}
\begin{proof}
  We have
  \begin{equation*}
    \epsilon_\lambda = \frac{\chi_\lambda(1)}{n!}\sum_{\mu\vdash n}\chi_\lambda(w_\mu)c_\mu.
  \end{equation*}
  By Theorem~\ref{theorem:basis-of-image}, the images of the functions $\epsilon_\lambda$ as $\lambda$ runs over $P(n,d)$ form a basis of the image of $K[S_d]$ in $\End_k V^{\otimes d}$.
  Replacing the vectors in a basis by non-zero scalar multiples still results in a basis, therefore the images of the irreducible characters
  \begin{equation}
    \label{eq:basis-of-chars}
    \{\chi_\lambda\mid \lambda \in P(n,d)\}
  \end{equation}
  also form a basis for the image of $K[S_d]$ in $\End_K V^{\otimes d}$.
  The transition matrix for going from the class functions $c_\lambda$ to the irreducible characters is the character table of $S_n$:
  \begin{equation*}
    \chi_\lambda = \sum_{\mu\vdash d} X_{\mu\lambda}c_\mu,
  \end{equation*}
  where $X_{\mu\lambda} = \chi_\lambda(w_\mu)$ for all partitions $\lambda$ and $\mu$ of $d$.
  Writing $Y = (Y_{\lambda\mu})$ for the inverse of the matrix $X$, we have:
  \begin{equation*}
    c_\lambda = \sum_{\mu\vdash d} Y_{\mu\lambda}\chi_\mu.
  \end{equation*}
  Taking images in $\End_KV^{\otimes d}$ of both sides of the above equation and using the second part of Theorem~\ref{theorem:basis-of-image} gives
  \begin{equation*}
    \bar c_\lambda = \sum_{\mu\in P(n,d)} Y_{\mu\lambda} \bar \chi_\mu \text{ for each } \lambda\in P(n,d).
  \end{equation*}
  The first part of Theorem~\ref{theorem:basis-of-image} now implies that 
  \begin{equation*}
    \{ \bar c_\lambda \mid \lambda\in P(n,d)\}
  \end{equation*}
  is a basis for the image of $K[S_d]$ in $\End_K V^{\otimes d}$ if and only if the $\# P(n,d)\times \# P(n,d)$ matrix $(Y_{\mu\lambda})_{\mu,\lambda\in P(n,d)}$ is non-singular.
  This is precisely Theorem~\ref{lemma:transmat} below.
\end{proof}
Besides being an essential ingredient in the proof of Theorem~\ref{theorem:class-basis-of-center}, the following theorem is interesting in its own right:
\begin{theorem}
  \label{lemma:transmat}
  Let $X_{\mu\lambda} = \chi_\lambda(w_\mu)$, the value of the character of the Specht module $V_\lambda$ at a permutation with cycle type $\mu$.
  Consider the matrix $X = (X_{\mu\lambda})$ whose rows and columns are indexed by partitions of $d$, and let $Y = X^{-1}$.
  Then the submatrix of $Y$ obtained by choosing only rows and columns whose partitions have at most $n$ parts is non-singular.
\end{theorem}
\begin{proof}
  Denote the power-sum symmetric functions and the Schur function corresponding to the partition $\lambda$ by $p_\lambda$ and $s_\lambda$ respectively.
  By the Murnaghan-Nakayama rule \cite[Eq.~(5.22)]{rtcv},
  \begin{equation}
    \label{eq:transition}
    p_\mu = \sum_{\lambda\vdash d} X_{\mu\lambda}s_\lambda.
  \end{equation}
  For a homogeneous symmetric function $f$ (which may be a priori viewed as a formal polynomial in infinitely many variables, or at least more variables than its degree), let $f(x_1,\dotsc,x_n)$ denote the specialization to $n$ variables.
  It is well known that $s_\lambda(x_1,\dotsc,x_n)=0$ if $\lambda$ has more than $n$ parts, and that the set
  \begin{equation*}
    \{s_\lambda(x_1,\dotsc,x_m)\mid \lambda\in P(n,d)\}
  \end{equation*}
  forms a basis of the space of symmetric functions of degree $d$ in $n$ variables \cite[Theorem~5.4.3]{rtcv}.
  Upon specialization, \eqref{eq:transition} becomes
  \begin{equation*}
    p_\mu(x_1,\dotsc,x_n) = \sum_{\lambda\in P(n,d)} X_{\mu\lambda}s_\lambda(x_1,\dotsc,x_n).
  \end{equation*}
  Since the specializations
  \begin{equation*}
    \{p_\lambda(x_1,\dotsc,x_n)\mid \lambda\in P(n,d)\}
  \end{equation*}
  of power sum symmetric functions also form a basis for symmetric functions of degree $d$ in $n$ variables \cite[Theorem~5.3.9]{rtcv}, it follows  that the submatrix of $X$ obtained by choosing the rows and columns corresponding to partitions in $P(n,d)$ is non-singular.

  By Schur's orthogonality relations and the fact that every permutation is conjugate to its own inverse, it follows that the character table matrix $X$ of $S_n$ satisfies $X'ZX = I$, where $Z$ is the diagonal matrix whose diagonal entry corresponding to $\lambda$ is the cardinality of the centralizer of a permutation with cycle type $\lambda$.
  It follows that the inverse $Y$ of $X$ is $X'Z^{-1}$.
  Thus, if the submatrix of $X$ obtained by selecting rows and columns indexed by partitions in $P(n,d)$ is non-singular, then so is the corresponding submatrix of $Y$.
\end{proof}
\begin{example}
  Consider the case where $n=2$ and $d = 3$.
  For each weak composition $\lambda$ the matrices $D$ with row and column sums $\lambda$, and the corresponding pairs $(\ii, \jj)$ with $D(\ii,\jj) = D$ are given in the following table:
  \begin{equation*}
    \begin{array}{ccc}
      \hline
      \lambda & \text{Possibilities for $D$} & (\ii, \jj)\\
      \hline
      (3,0) & \mat 3000 & ((1, 1, 1), (1, 1, 1))\\[0.5cm]
      (0,3) & \mat 0003 & ((2, 2, 2), (2, 2, 2))\\[0.5cm]
      (2,1) & \mat 1110, \mat 2001 & ((1, 1, 2), (1, 2, 1)), ((1, 1, 2), (1, 1, 2))\\[0.5cm]
      (1,2) & \mat 0111, \mat 1002 & ((1, 2, 2), (2, 1, 2)), ((1, 2, 2), (1, 2, 2))\\[0.5cm]
      \hline
    \end{array}
  \end{equation*}
  The partitions in $P(2,3)$ are $(3)$ and $(2,1)$.
  Thus the centre of $S_K(2,3)$, has basis elements $\bar c_{(3)}$ and $\bar c_{(2,1)}$ whose expansions in terms of Schur's basis $\xi_D$ are given by:
  \begin{align*}
    \bar c_{(3)} & = 2\xi_{\smat 3000} + 2\xi_{\smat 0003} + \xi_{\smat 1110} + \xi_{\smat 0110}\\[0.3cm]
    \bar c_{(2, 1)} & = 3\xi_{\smat 3000} + 3\xi_{\smat 0003} + \xi_{\smat 1110} + \xi_{\smat 2001} + \xi_{\smat 0111} + \xi_{\smat 1002}\\[0.3cm]
  \end{align*}
  We also have
  \begin{equation*}
    \bar c_{(1,1,1)} = \xi_{\smat 3000} + \xi_{\smat 0003} + \xi_{\smat 2001} + \xi_{\smat 1002} = \bar c_{(2,1)} - \bar c_{(3)}.
  \end{equation*}
  And, finally (using character values for $S_3$) the primitive central idempotents of $S_K(n,d)$ are given by:
  \begin{align*}
    \bar \epsilon_{(3)} & = \bar c_{(3)} + \bar c_{(2,1)} + \bar c_{(1,1,1)} = 2\bar c_{(2,1)}\\
    \bar \epsilon_{(2,1)} & = -\bar c_{(3)} + 2 \bar c_{(1,1,1)} = \bar c_{(2,1)} + \bar c_{(3)}.
  \end{align*}
  Also, we find (as expected) that
  \begin{equation*}
    \bar \epsilon_{(1,1,1)}  = \bar c_{(3)} - \bar c_{(2,1)} + \bar c_{(1,1,1)} = 0 
  \end{equation*}
\end{example}
\section*{Acknowledgements}
The authors are very grateful to Vijay Kodiyalam and Steffen K\"onig for their help and encouragement. Geetha was supported by the Institute of Mathematical Sciences, Chennai, while working on this article.

\bibliographystyle{plain}
\bibliography{center}

\begin{thebibliography}{1}

\bibitem{structure}
T.~Geetha and Amritanshu Prasad.
\newblock Graphic interpretation of the structure constants of the {S}chur
  algebra.
\newblock In {\em Electronic booklet of Proceedings of the International
  Congress of Women Mathematicians}, 2014.
\newblock Abstract no. 20140020; also available at
  \url{http://arxiv.org/abs/1409.1414}.

\bibitem{Green}
James~A. Green.
\newblock {\em Polynomial Representations of $GL_n$}, volume 830 of {\em
  Lecture Notes in Mathematics}.
\newblock Springer-Verlag, Berlin Heidelberg, 1980.

\bibitem{rtcv}
Amritanshu Prasad.
\newblock {\em Representation Theory: A Combinatorial Viewpoint}.
\newblock Cambridge Studies in Advanced Mathematics, 147. Cambridge University
  Press, 2014.

\end{thebibliography}
\end{document}